\documentclass[leqno,12pt]{amsart}
\usepackage{latexsym}
\usepackage{mathrsfs}
\usepackage{amsmath}
\usepackage{amssymb}
\usepackage[bookmarks]{hyperref}
\usepackage{pstricks,pst-plot}

\font\tenscr=rsfs10 
\font\sevenscr=rsfs7 
\font\fivescr=rsfs5 
\skewchar\tenscr='177 \skewchar\sevenscr='177 \skewchar\fivescr='177
\newfam\scrfam \textfont\scrfam=\tenscr \scriptfont\scrfam=\sevenscr
\scriptscriptfont\scrfam=\fivescr
\def\scr{\fam\scrfam}

\font\tenscr=rsfs10 
\font\sevenscr=rsfs7 
\font\fivescr=rsfs5 
\skewchar\tenscr='177 \skewchar\sevenscr='177 \skewchar\fivescr='177
\newfam\scrfam \textfont\scrfam=\tenscr \scriptfont\scrfam=\sevenscr
\scriptscriptfont\scrfam=\fivescr
\def\scr{\fam\scrfam}

\newtheorem{theorem}{Theorem}[section]

\newtheorem{corollary}[theorem]{Corollory}
\newtheorem{lemma}[theorem]{Lemma}

\newtheorem{definition} [theorem]{Definition}

\newcommand{\bthm}{\begin{theorem}}
\newcommand{\ethm}{\end{theorem}}
\newcommand{\blem}{\begin{lemma}}
\newcommand{\elem}{\end{lemma}}
\newcommand{\bcor}{\begin{corollary}}
\newcommand{\ecor}{\end{corollary}}
\newcommand{\bprop}{\begin{proposition}}
\newcommand{\eprop}{\end{proposition}}
\newcommand{\bdefn}{\begin{definition}}
\newcommand{\edefn}{\end{definition}}
\newcommand{\bpf}{\begin{proof}}
\newcommand{\epf}{\end{proof}}

\newcommand{\bi}{\begin{itemize}}
\newcommand{\ei}{\end{itemize}}

\newcommand{\bc}{\begin{cases}}
\newcommand{\ec}{\end{cases}}

\newcommand{\ba}{\begin{array}}
\newcommand{\ea}{\end{array}}

\newcommand{\be}{\begin{equation}}
\newcommand{\ee}{\end{equation}}

\newcommand{\bea}{\begin{eqnarray}}
\newcommand{\eea}{\end{eqnarray}}

\newcommand{\beaa}{\begin{eqnarray*}}
\newcommand{\eeaa}{\end{eqnarray*}}

\newcommand{\beastar}{\begin{eqnarray*}}
\newcommand{\eeastar}{\end{eqnarray*}}

\begin{document}

\def\tA{\tilde A}
\def\tX{\tilde X}
\def\tf{\tilde f}
\def\tpi{\tilde \pi}
\def\th{\tilde h}
\def\ta{\tilde \alpha}

\def \vep {\varepsilon}

\def \cd {, \ldots ,}
\def \lf{\| f \|}
\def \bs{\setminus}
\def \ep{\varepsilon}
\def \sig{\Sigma}
\def \si {\sigma}
\def \gam {\gamma}
\def \cinf {C^\infty}
\def \cid {C^\infty (\partial D)}

\def\h#1{\widehat {#1}}
\def\hh#1{\widehat {#1} \bs  {#1}}
\def\hk#1#2{{\widehat {#2}}^{#1}}
\def\hhk#1#2{{\widehat {#2}}^{#1} \bs {#2}}

\def\hr#1{h_r({#1})}
\def\hhr#1{h_r({#1}) \bs {#1}}
\def\hrk#1#2{h_r^{#1}({#2})}
\def\hhrk#1#2{h_r^{#1}({#2}) \bs {#2}}

\def\kh#1#2{{}^{#1}{\widehat {#2}}}
\def\khr#1#2{{}^{#1}h_r({#2})}
\def\kshr#1#2{{}_{#1}h_r({#2})}

\def\<{\langle}
\def\>{\rangle}

\def\spshell{\{ z\in \CN: 1<\|z\| < \rho\}}

\def \C {\mathbb C}
\def \CN {\mathbb C^N}
\def\R {\mathbb R}

\def \Z {\mathbb Z}

\def \bbr {\mathbb R}
\def \bbrs {\mathbb R^2}
\def \bbrn {\mathbb R^n}

\def \cc {\mathscr C}
\def \kk {\mathscr K}
\def \mm {\mathfrak M}

\def \od {\overline D}
\def \ol {\overline L}
\def \oj {\overline J}

\def\sC{{\scr C}}
\def\sF{{\scr F}}
\def\sG{{\scr G}}

\def \cn {c_1, \ldots , c_n}
\def \zx { z \in \widehat X }
\def \zxx {\{ z \in \widehat X \} }
\def \siv {\sum^\infty _{j=1}}
\def \sivv {\sum^\infty _{j=n+1}}
\def \smv {\sum^m_{j=n+1}}
\def \snv {\sum^n_{j=j_0+1}}
\def \bjkz {B_{jk} (z_0)}
\def \epp {\varepsilon^2_n/2}
\def \epn {\varepsilon_{n+1}}
\def \nn { {n+1} }
\def \nnn {{2N}}
\def \ff {F_{z_{0}}}
\def \fk {f_{k_{1}}}
\def \fkv {f_{k_{v}}}
\def \tf {\tilde{f}}

\def \ma {\mathfrak{M}_A}  
\def \mb {\mathfrak{M}_B}
\def \muc {\mathfrak{M}_{\mathscr U}}
\def \uc {\mathscr U}
\def \xy {(x,y)}
\def \fg {f \otimes g}
\def \vp {\varphi}
\def \pl {\partial L}
\def \ad {A(D)}
\def \aid {A^\infty (D)}
\def \ot {\otimes}
\def \rbl {R_b (L)}
\def\Int {{\rm Int}}

\def\rowonly#1#2{#1_1,\ldots,#1_#2}
\def\row#1#2{(#1,\ldots,#1_#2)}
\def\diam{\mathop{\rm diam}\nolimits}

\def\pB{\partial B}
\def\oB{\overline B}

\def\endhat{\widehat{\phantom j}}
\def\endhatk{\widehat{\phantom j}{\phantom |}^k}

\def\pji{p_j^{-1}}
\def\oh{\overline H}

\subjclass[2000]{32E20, 46J10, 46J15}
\title[a one-point part off the Shilov boundary]{A doubly generated uniform algebra\\ with a one-point Gleason part\\ off its Shilov boundary}
\author{Alexander J. Izzo}
\address{Department of Mathematics and Statistics, Bowling Green State University, Bowling Green, OH 43403}
\email{aizzo@bgsu.edu}

\begin{abstract}
It is shown that there exists a compact set $X$ in $\C^2$ with a nontrivial polynomial hull $\h X$ such that some point of $\hh X$ is a one-point Gleason part for $P(X)$.  Furthermore, $X$ can chosen so that $P(X)$ has a dense set of invertible elements.
\end{abstract}

\maketitle

 \section{Introduction}
 
Examples of uniform algebras $A$ with a one-point Gleason part lying off the Shilov boundary for $A$ have long been known.  See for instance \cite[p.~187]{Stout-old} for a well-known example defined on the 2-torus and known as the big
disc algebra.  Brian Cole's method of root extensions \cite{Cole} yields uniform algebras with nontrivial Shilov boundary such that \emph{every} point is a one-point Gleason part.  Recently, Cole, Ghosh, and the present author \cite{CGI} gave a triply generated example, that is, they constructed a compact set $X$ in $\C^3$ with a nontrivial polynomial hull $\h X$ such that every point of $\h X$ is a one-point Gleason part for $P(X)$.  This suggests the question of whether for $X$ a compact set in $\C^2$, the set $\hh X$ can contain a one-point Gleason part.  The main purpose of the present paper is to show that this can indeed happen.  Furthermore, $X$ can be chosen so that $P(X)$ has a dense set of invertible elements, and hence in particular, so that $\h X$ contains no analytic discs.  (As noted by Garth~Dales and Joel~Feinstein in \cite{DalesF}, the condition that $P(X)$ has a dense set of invertible elements is strictly stronger than the condition that $\h X$ contains no analytic disc.)  Since for $X$ a compact set in $\C^1$, there can never be a one-point Gleason part for $P(X)$ lying in $\hh X$, this result is, in a certain sense, optimal.  Whether the result can be strengthened to obtain a compact set $X$ in $\C^2$ with nontrivial polynomial hull such that \emph{every} point of $\h X$ is a one-point Gleason part for $P(X)$ remains open.

We will denote the open unit ball in $\CN$ by $B$ and the boundary of an open set $\Omega$ in $\CN$ by $\partial\Omega$.
 
 \bthm\label{one-point part}
 There exists a compact set $X\subset\pB\subset\C^2$ such that the origin is in $\hh X$ and is a one-point Gleason part for $P(X)$. 
 \ethm
 
Following Dales and Feinstein, 
we will say that a uniform algebra $A$ has dense invertibles if the invertible elements of $A$ are dense in $A$.
As alluded to earlier, the set $X$ in the above theorem can be chosen so that $P(X)$ has a dense invertibles.  
Using extensions of the notions of polynomial and rational hulls introduced by the author in \cite{Izzo-spaces}, we will establish the following more general result.  The definitions of these new notions, the $k$-polynomial hull $\hk kX$ of $X$ and the $k$-rational hull $\hrk kX$ of $X$, will be recalled below in the next section.
 
 \bthm\label{dense invert}
 Let $\Sigma\subset\CN$ be any compact set with nontrivial $k$-polynomial hull with $k\geq 2$, and let $x_0\in \hk k\Sigma$.  Then $\Sigma$ contains a compact set $X$ such that $x_0$ is in $\hhrk {k-1}X\subset \hh X$ and is a one-point Gleason part for $P(X)$, and $P(X)$ has dense invertibles. 
 \ethm
 
 The following corollary is an almost immediate consequence.
 
 \bcor\label{boundary}
 Let $\Omega\subset \CN$ {\rm ($N\geq 2$)} be any bounded open set, and let $x_0\in \Omega$.  Then there exists a compact set $X\subset\partial \Omega$ such that $x_0$ is in $\hh X$ and is a one-point Gleason part for $P(X)$, and $P(X)$ has dense invertibles. 
 \ecor
 
Using results from the author's paper \cite{Izzo-spaces}, we will obtain the following as another corollary.
 
 \bcor\label{Cantor set}
 There exists a Cantor set $X$ in $\C^3$ such that $\hh X$ is nonempty, some point of $\hh X$ is a one-point Gleason part for $P(X)$, and $P(X)$ has dense invertibles. 
 \ecor

The existence of a Cantor set $X$ in $\C^3$ such that $\hh X$ is nonempty and $P(X)$ has dense invertibles was proved by the author as 
\cite[Theorem~1.4]{Izzo-spaces}.  What is new in Corollary~\ref{Cantor set} is that the Cantor set $X$ can be chosen so that some point of $\hh X$ is a one-point Gleason part for $P(X)$.  Whether there exists a Cantor set $X$ in $\C^3$ with nontrivial polynomial hull so that \emph{every} point of $\h X$ is a one-point Gleason part remains open.  However, the author showed in \cite[Theorem~7.2]{Izzo-spaces} that a Cantor set with these properties does exist in $\C^4$.

In connection with the above results, we mention that, as noted by Dales and Feinstein~\cite{DalesF} (and implicitly noted by Stolzenberg~\cite{Stol1}), the polynomial and rational hulls of a compact set $X$ coincide whenever $P(X)$ has dense invertibles.

In the next section, in addition to giving the definitions of the $k$-polynomial and $k$-rational hulls, we make explicit some standard definitions and notations already used above, and we recall some known results we will need.  The proofs of the results stated above are given in Section~3.

%

\section{Preliminaries}~\label{prelim}

For
$X$ a compact Hausdorff space, we denote by $C(X)$ the algebra of all continuous complex-valued functions on $X$ with the supremum norm
$ \|f\|_{X} = \sup\{ |f(x)| : x \in X \}$.  A \emph{uniform algebra} on $X$ is a closed subalgebra of $C(X)$ that contains the constant functions and separates
the points of $X$.  

For a compact set $X$ in $\CN$, the \emph{polynomial hull} $\h X$ of $X$ is defined by
$$\h X=\{z\in\CN:|p(z)|\leq \max_{x\in X}|p(x)|\
\mbox{\rm{for\ all\ polynomials}}\ p\},$$
and the
\emph{rational hull} $\hr X$ of $X$ is defined by
$$\hr X = \{z\in\C^N: p(z)\in p(X)\ 
\mbox{\rm{for\ all\ polynomials}}\ p
\}.$$
An equivalent formulation of the definition of $\hr X$ is that $\hr X$ consists precisely of those points $z\in \C^N$ such that every polynomial that vanishes at $z$ also has a zero on $X$.

Observe that when $N=2$, the statement $z\in \hr X$ means that every analytic subvariety of $\C^N$ (of pure positive dimension) passing through $z$ intersects $X$, whereas when $N>2$, the statement $z\in \hr X$ means only that every pure codimension 1 analytic subvariety of $\C^N$ passing through $z$ intersects $X$.  This observation suggests the following extensions of the notions of polynomial and rational hulls introduced by the author in \cite[Definitions~3.1 and~3.2]{Izzo-spaces}.  (Here and throughout the paper, by an analytic subvariety of $\CN$ of pure codimension $m$, we mean a pure dimensional analytic subvariety of $\CN$ of pure dimension $N-m$.)
For $1\leq k\leq N$, the \emph{$k$-rational hull} $\hrk kX$ of a compact set
$X\subset\CN$ is defined by
\begin{eqnarray*}
\hrk kX&\!\!\!=&\!\!\!\{ z\in \CN: {\rm every\ analytic\ subvariety\ of\ } \CN {\rm\ of\ pure\ codimension}\leq k \\
& &\qquad\qquad  {\rm that\ passes\ through\ } z\ {\rm intersects\ } X\}.
\end{eqnarray*}
For $2\leq k\leq N$, the \emph{$k$-polynomial hull} $\hk kX$ of 
$X\subset\CN$ is defined by
\begin{eqnarray*}
\hk kX&\!\!\!=&\!\!\!\{ z\in \C^N: z\in \hrk {k-1}X\ {\rm and\ } z\in \h {X\cap V}\ {\rm for\ every} \\ 
& &\qquad\qquad{\rm analytic\ subvariety}\ V\ {\rm of\ } \CN {\rm\ of\ pure\  codimension}\leq k-1 \\
& &\qquad\qquad  {\rm that\ passes\ through\ } z\} .
\end{eqnarray*}
The \emph{$1$-polynomial hull} $\hk 1X$ of $X$ is defined to be the usual\vadjust{\kern 2pt} polynomial hull $\h X$.

It is immediate from the definitions that 
$$\h X=\hk 1X\supset \hr X=\hrk 1X\supset \hk 2X\supset \hrk 2X 
\supset\cdots\supset \hk NX \supset \hrk NX=X.$$
The set $X$ is said to be \emph{$k$-polynomially convex} if $\hk kX = X$ and \emph{$k$-rationally convex} if $\hrk kX = X$.
It is easily verified that the $k$-polynomial hull of a compact set is $k$-polynomially convex and the $k$-rational hull is $k$-rationally convex.

One could also consider modifications of the above definitions.
For $1\leq k\leq N$, the \emph{quasi-$k$-rational hull} $\khr kX$ of 
$X$ is defined by
\begin{eqnarray*}
\khr kX&\!\!\!=&\!\!\!\{ z\in \CN: {\rm if}\ p_1,\ldots, p_k\ {\rm are\ polynomials\ such\ that}  \\
&&p_1(z)=0,\ldots,  p_k(z)=0,\ {\rm then}\ p_1,\ldots, p_k\ {\rm have\ a\ common\ zero\ on}\ X\}.
\end{eqnarray*}
For $2\leq k\leq N$, the \emph{quasi-$k$-polynomial hull} $\kh kX$ of 
$X$ is defined by
\begin{eqnarray*}
\kh kX&\!\!\!=&\!\!\!\ \bigl\{ z\in \C^N: z\in \khr {k-1}X\ {\rm and\ } {\rm if}\ p_1,\ldots, p_{k-1}\ {\rm are\ polynomials\ such\ that}  \\
&&p_1(z)=0,\ldots, p_{k-1}(z)=0,\ {\rm then}\ z\in (X \cap\{p_1=0, \ldots,  p_{k-1}=0\})\endhat\,\bigr\}.
\end{eqnarray*}
The \emph{quasi-$1$-polynomial hull} $\kh 1X$ of $X$ is defined to be the usual polynomial hull $\h X$.

Note that the $k$-hulls are ostensibly smaller than the corresponding quasi-$k$-hulls.  We will use exclusively the $k$-hulls.
However, the quasi-$k$-hulls are perhaps more intuitive than the $k$-hulls, and the reader who wishes to do so, may replace the $k$-hulls by the quasi-$k$-hulls throughout.  

We say that a hull $X_H$ is {\it nontrivial} if the set $X_H\setminus X$ is nonempty.

We denote by 
$P(X)$ the uniform closure on $X\subset\CN$ of the polynomials in the complex coordinate functions $z_1,\ldots, z_N$, and we denote by $R(X)$ the uniform closure of the rational functions  holomorphic on (a neighborhood of) $X$. 
Both $P(X)$ and $R(X)$ are uniform algebras, and
it is well known that the maximal ideal space of $P(X)$ can be naturally identified with $\h X$, and the maximal ideal space of $R(X)$ can be naturally identified with $\hr X$.

By an \emph{analytic disc} in $\CN$, we mean an injective holomorphic map $\sigma: \{ z\in \C: |z|< 1\}\rightarrow\CN$.
By the statement that a subset $S$ of $\CN$ contains no analytic discs, we mean that there is no analytic disc in $\CN$ whose image is contained in $S$.

Let $A$ be a uniform algebra on a compact space $X$.
The \emph{Gleason parts} for the uniform algebra $A$ are the equivalence classes in the maximal ideal space of $A$ under the equivalence relation $\varphi\sim\psi$ if $\|\varphi-\psi\|<2$ in the norm on the dual space $A^*$.  (That this really is an equivalence relation is well-known but {\it not\/} obvious!)
An alternative formulation of the definition of Gleason part is that two points $\varphi$ and $\psi$ lie in the same Gleason part if and only if there exists a constant $c>0$ such that 
\[ 1/c < u(\varphi)/u(\psi) < c \]
for every function $u>0$ that is the real part of a function in $A$.  See 
\cite[Theorem~VI.2.1]{Gamelin} for the equivalence of these two formulations.  It is easily seen that in the second formulation, the condition $u>0$ can be relaced by $u\geq 0$, and it is obvious that $u\geq 0$ can be replaced by $u\leq 0$.
We say that a Gleason part is \emph{nontrivial} if it contains more than one point.
It is immediate that the presence of an analytic disc in the maximal ideal space of $A$ implies the existence of a nontrivial Gleason part.

The real part of a complex number (or function) $z$ will be denoted by $\Re z$.   

For the reader's convenience we recall here some lemmas that we will use.  The first of these is standard, and a short proof can be found in \cite[Lemma~1.7.4]{Stout-new}.  The others are proven in \cite{Izzo-spaces}.

\begin{lemma}\label{Stout-lemma}
If $X\subset \CN$ is a polynomially convex set, and if $E\subset X$ is polynomially convex, then for every holomorphic function $f$ defined on a neighborhood of $X$, the set $E\cup (X\cap f^{-1}(0))$ is polynomially convex.
\end{lemma}

\begin{lemma} \label{hullintersection} \cite[Lemma~4.1]{Izzo-spaces}
Let $\kk$ be a collection of compact sets in $\CN$ totally ordered by inclusion.  Let $K_\infty=\bigcap_{K\in \kk} K$.  Then $\hk k {K}_\infty=\bigcap_{K\in \kk} \hk k K$ and $\hrk k {K_\infty} = \bigcap_{K\in \kk} \hrk k K$.
\end{lemma}

Recall that a subset of a space is called \emph{perfect} if it is closed and has no isolated points.  Every space contains a unique largest perfect subset (which can be empty), namely the closure of the union of all perfect subsets of the space.

\begin{lemma}\label{reductiontop} \cite[Lemma~4.2]{Izzo-spaces}
Let $X\subset \CN$ be a 
compact set and let $E$ be the largest perfect subset of $X$.  Then $\hhk k E \supset \hhk k X$ and $\hhrk k E\supset\hhrk k X$.
\end{lemma}

The next lemma plays a key role in our proofs.

\begin{lemma} \label{foundation} \cite[Lemma~5.9]{Izzo-spaces}
Let $\Sigma\subset\C^N$ be a compact set, let $p$ be a polynomial on $\C^N$, and let $X=\{\Re p\leq 0\} \cap \Sigma$.  Let $k\geq 2$ be an integer.  Then 
$\{\Re p\leq 0\} \cap \hk {k-1}\Sigma \supset \hk {k-1}X \supset  \hrk {k-1}X \supset
\{\Re p\leq 0\} \cap \hk k\Sigma$.
\end{lemma}

Our final lemma is a special case of the previous one because, as is shown in the proof of Corollary~\ref{boundary} below,
$\h {\partial B} = \hk N{\partial B}=\oB$.  The reader who wishes can, however, supply a direct proof of the lemma similar to the proof of \cite[Lemma~3.2]{Izzo}.

\begin{lemma}\label{hullofcap} \cite[Lemma~6.3]{Izzo-spaces}
Let $N\geq 2$.  Let $p$ be a polynomial on $\CN$, and let $X=\{\Re p\leq 0\} \cap \pB$.  Then $\h X = \hrk {N-1} X = \{\Re p\leq 0\} \cap \oB$.
\end{lemma}

%

\section {The proofs}

In this section we prove the results stated in the introduction. 
 
Theorem~\ref{one-point part} is of course contained in Corollary~\ref{boundary}.  However, we will give a direct proof so as to present the construction giving the one-point Gleason part in its simplest form without the additional complications involved in obtaining dense invertibles.

The proof of Theorem~\ref{one-point part} has some similarities with the construction of a polynomial hull without analytic discs given by Julien~Duval and Norman~Levenberg~\cite{DuvalL} (or see \cite[Lemma~1.7.5]{Stout-new}).  The proof of Theorem~\ref{dense invert} involves combining the ingredients in the proof of Theorem~\ref{one-point part} with ingredients in the proof of \cite[Theorem~5.2]{Izzo-spaces}.  (Theorem~5.2 of
\cite{Izzo-spaces} generalizes \cite[Theorem~4.1]{Izzo} whose proof was inspired by the constructions of Duval and Levenberg~\cite{DuvalL} and of Dales and Feinstein \cite{DalesF}.)
 
\begin{proof}[Proof of Theorem \ref{one-point part}]
Choose a sequence $\{E_j\}_{j=1}^\infty$ of compact polynomially convex subsets of $\oB\setminus\{0\}$ such that every point of $\oB\setminus\{0\}$ lies in $E_j$ for infinitely many values of $j$.  We will construct a sequence of polynomials $\{f_j\}_{j=1}^\infty$ such that the sets
\be
X_j=\{\Re f_j\leq 0\}\cap \pB\label{Xj}
\ee
form a decreasing sequence and such that for each $j$ we have 
\[ \hbox{$f_j(0)=-1$\quad and\quad $\Re f_j>-1/j$ on $E_j\cap {\h X}_{j+1}$}.\]
Then letting $X=\bigcap_{j=1}^\infty X_j$ we have $\h X=\bigcap_{j=1}^\infty {\h X}_j=\bigcap_{j=1}^\infty \bigl(\{\Re f_j\leq 0\}\cap \oB\bigr)$ by Lemmas~\ref{hullintersection} and~\ref{hullofcap}.  In particular, $0\in\h X$.  Furthermore, since each point of $\h X\setminus\{0\}$ lies in $E_j$ for infinitely many $j$, for each point $x\in \h X\setminus \{0\}$ there are infinitely many $j$ such that $\Re f_j(x)> -1/j$, and hence such that $0\leq\Re f_j(x) / \Re f_j(0) < 1/j$.  Since $\Re f_j\leq 0$ on $\h X$ for each $j$, this gives that $\{0\}$ is a one-point Gleason part for $P(X)$ by the remarks in Section~2.

We construct the $f_j$ inductively.  To begin, set $f_1$ identically equal to $-1$.  Then assume for the purpose of induction that polynomials $f_1,\ldots, f_n$ have been chosen such that the sets $X_1,\ldots, X_n$ defined as in (\ref{Xj}) for $j=1,\ldots, n$ form a decreasing sequence, for each $j=1,\ldots, n$ we have $f_j(0)= -1$, and for each $j=1,\ldots, n-1$ we have $\Re f_j> -1/j$ on $E_j\cap {\h X}_{j+1}$.  Let $L_n=\{\Re f_n\geq 0\} \cap \oB$ and $C_n=E_n\cap \{\Re f_n\leq -1/n\}$.  The sets $L_n$ and $C_n$ are each polynomially convex.  Because the sets $f_n(L_n)$ and $f_n(C_n)$ lie in disjoint half-planes, their polynomial hulls are disjoint.  Therefore, $L_n \cup C_n$ is polynomially convex by Kallin's lemma \cite{Kallin} (or see \cite[Theorem~1.6.19]{Stout-new}).
Adjoining a single point to a polynomially convex set yields another polynomially convex set \cite[Lemma~2.2]{AIW2001}, so $L_n\cup C_n \cup \{0\}$ is also polynomially convex.  Therefore, there exists a polynomial $f_{n+1}$ such that $f_{n+1}(0)= -1$ and $\Re f_{n+1} >0$ on $L_n\cup C_n$.  Now define $X_{n+1}$ as in (\ref{Xj}) with $j=n+1$ and observe that then $X_n\supset X_{n+1}$ and $\Re f_n > -1/n$ on $E_n\cap {\h X}_{n+1}$, so the induction can continue.
\epf

\bpf[Proof of Theorem \ref{dense invert}]
Choose a sequence $\{E_j\}_{j=1}^\infty$ of compact polynomially convex subsets of $\h \Sigma \setminus \{x_0\}$ such that every point of $\h\Sigma\setminus\{x_0\}$ lies in $E_j$ for infinitely many values of $j$.  
Also choose a sequence of polynomials $\{p_j\}_{j=1}^\infty$ that is dense in $P(\Sigma)$ and such that $p_j(x_0)\neq 0$ for each $j$.
Let $Z_j=\h \Sigma \cap p_j^{-1}(0)$ for each $j$.
We will construct a sequence of polynomials $\{f_j\}_{j=1}^\infty$ such that the sets
\be
X_j=\{\Re f_j\leq 0\}\cap \Sigma\label{XSigma}
\ee
form a decreasing sequence and such that for each $j$ we have 
\[ \hbox{$f_j(x_0)=-1$,\quad $\Re f_j>0$ on $Z_j$,\quad and\quad $\Re f_j>-1/j$ on $E_j\cap {\h X}_{j+1}$}.\]
Then letting $X=\bigcap_{j=1}^\infty X_j$, we have by Lemmas~\ref{hullintersection} and~\ref{foundation} that $\hrk {k-1}X=\bigcap_{j=1}^\infty {\hrk {k-1}{X_j}}\supset\bigcap_{j=1}^\infty \bigl(\{\Re f_j\leq 0\}\cap \hk k\Sigma\bigr)$.  In particular, $x_0\in\hrk {k-1}X$.  
Also $\h X=\bigcap_{j=1}^\infty {\h X}_j\subset \bigcap_{j=1}^\infty \bigl(\{\Re f_j\leq 0\}\cap \h\Sigma\bigr)$, again by 
Lemmas~\ref{hullintersection} and~\ref{foundation}.  Thus $\h X$ is disjoint from each $Z_j$.
Consequently, each $p_j$ is invertible in $P(X)$, and hence $P(X)$ has dense invertibles.
Furthermore, since each point of $\h \Sigma\setminus\{x_0\}$ lies in $E_j$ for infinitely many $j$, for each point $x\in \h \Sigma\setminus \{x_0\}$ there are infinitely many $j$ such that $\Re f_j(x)> -1/j$, and hence such that $0\leq\Re f_j(x) / \Re f_j(x_0) < 1/j$.  Since $\Re f_j\leq 0$ on $\h X$ for each $j$, this gives that $\{x_0\}$ is a one-point Gleason part for $P(X)$ by the remarks in Section~2. 

We construct the $f_j$ inductively.  The set $Z_1$ is polynomially convex, and hence $Z_1\cup\{x_0\}$ is also polynomially convex (by \cite[Lemma~2.2]{AIW2001}), so there is a polynomial $f_1$ such that
\[ \hbox{$f_1(x_0)=-1$\quad and \quad$\Re f_1>0$ on $Z_1$.}\]
Set $X_1=\{\Re f_1\leq 0\}\cap \Sigma$.  For the inductive step, assume that polynomials $f_1,\ldots, f_n$ have been chosen such that the sets $X_1,\ldots, X_n$ defined as in (\ref{XSigma}) for $j=1,\ldots, n$ form a decreasing sequence, for each $j=1,\ldots, n$ we have $f_j(x_0)= -1$ and $\Re f_j>0$ on $Z_j$, and for each $j=1,\ldots, n-1$ we have $\Re f_j> -1/j$ on $E_j\cap {\h X}_{j+1}$.  Let $L_n=\{\Re f_n\geq 0\} \cap \h\Sigma$ and $C_n=E_n\cap \{\Re f_n\leq -1/n\}$.  The sets $L_n$ and $C_n$ are each polynomially convex.  Because the sets $f_n(L_n)$ and $f_n(C_n)$ lie in disjoint half-planes, their polynomial hulls are disjoint.  Therefore, $L_n \cup C_n$ is polynomially convex by Kallin's lemma \cite{Kallin} (or see \cite[Theorem~1.6.19]{Stout-new}).
Applying Lemma~\ref{Stout-lemma} now gives that $L_n \cup C_n \cup Z_{n+1}$ is polynomially convex.
Consequently, $L_n\cup C_n \cup Z_{n+1} \cup \{x_0\}$ is also polynomially convex (by  \cite[Lemma~2.2]{AIW2001}).  Therefore, there exists a polynomial $f_{n+1}$ such that $f_{n+1}(x_0)= -1$ and $\Re f_{n+1} >0$ on $L_n\cup C_n \cup Z_{n+1}$.  Now set $X_{n+1}=\{\Re f_{n+1}\leq 0\}\cap \Sigma$.
Because $\Re f_{n+1}>0$ on $L_n$, we have $X_n\supset X_{n+1}$, and because $\Re f_{n+1}>0$ on $C_n$, we have $\Re f_n> -1/n$ on $E_n\cap {\h X}_{n+1}$.
Thus the induction can continue.
\epf

\bpf[Proof of Corollary~\ref{boundary}]
Let $x\in \Omega$ be arbitrary.  Then every analytic subvariety $V$ of $\CN$ of pure positive dimension passing through $x$ intersects $\partial\Omega$.  Consequently, $x$ lies in $\hrk {N-1}{\partial\Omega}$, and applying the maximum principle to the irreducible component of $V$ through $x$ shows that $x$ is in $\h {\partial\Omega\cap V}$.  Thus $\Omega\subset \hk N{\partial\Omega}$.  The corollary now follows immediately from Theorem~\ref{dense invert}.
\epf

\bpf[Proof of Corollary~\ref{Cantor set}]
By \cite[Theorem~6.1]{Izzo-spaces}
there exists a Cantor set $\Sigma$ in $\C^3$ such that $\hk 2\Sigma$ contains the closed unit ball of $\C^3$.  Let $x_0$ be a point of $\hhk 2\Sigma$.  Theorem~\ref{dense invert} gives that $\Sigma$ contains a compact set $K$ such that $x_0$ is in $\hh K$ and is a one point Gleason part for $P(K)$, and $P(K)$ has dense invertibles.  Let $X$ be the largest perfect subset of $K$.  Then $\hh X \supset \hh K$ by Lemma~2.3, so $x_0$ is in $\hh X$, and the condition that $x_0$ is a one-point Gleason part for $P(K)$ implies that $x_0$ is also a 
one-point Gleason part for $P(X)$.  The set $X$ is a Cantor set since it is a perfect subset of the Cantor set $\Sigma$.  Finally density of the invertibles in $P(K)$ implies density of the invertibles in $P(X)$.
\epf

\end{document}